\documentclass[a4paper,12pt]{amsart}
\usepackage{amssymb}
\usepackage{ifthen}
\usepackage{graphicx}
\usepackage[usenames]{color}
\nonstopmode
\newtheorem{thm}{Theorem}
\newtheorem{cor}{Corollary}
\newtheorem{lem}{Lemma}

\newtheorem{claim}{Claim}

\newtheorem{conj}{Conjecture}
\newtheorem{prob}{Problem}

\theoremstyle{definition}
\newtheorem{defn}{Definition}
\newtheorem{example}{Example}

\newenvironment{rem}{%
\bigskip
\noindent \textsl{{\sl Remark. }}}{\bigskip}
\newenvironment{rems}{%
\bigskip
\noindent \textsl{{\sl Remarks. }}}{\bigskip}

\newenvironment{pf}[1][]{%
 \vskip 1mm
 \noindent
 \ifthenelse{\equal{#1}{}}%
  {{\slshape Proof. }}%
  {{\slshape #1.} }%
 }%
{\qed\medskip}

\newcounter{alphabet}
\newcounter{tmp}
\newenvironment{Thm}[1][]{\refstepcounter{alphabet}%
\bigskip%
\noindent%
{\bf Theorem \Alph{alphabet}}%
\ifthenelse{\equal{#1}{}}{}{ (#1)}%
{\bf .} \itshape}{\vskip 8pt}

\makeatletter
\newcommand{\Ref}[1]{\@ifundefined{r@#1}{}{\setcounter{tmp}{\ref{#1}}\Alph{tmp}}}
\makeatother

\newenvironment{Lem}[1][]{\refstepcounter{alphabet}%
\bigskip%
\noindent%
{\bf Lemma \Alph{alphabet}}%
{\bf .} \itshape}{\vskip 8pt}

\newcommand{\IR}{{\mathbb R}}

\newcommand{\IC}{{\mathbb C}}
\newcommand{\ID}{{\mathbb D}}

\def\be{\begin{equation}}
\def\ee{\end{equation}}

\newcommand{\bee}{\begin{enumerate}}
\newcommand{\eee}{\end{enumerate}}

\newcommand{\blem}{\begin{lem}}
\newcommand{\elem}{\end{lem}}
\newcommand{\bthm}{\begin{thm}}
\newcommand{\ethm}{\end{thm}}
\newcommand{\bcor}{\begin{cor}}
\newcommand{\ecor}{\end{cor}}
\newcommand{\beg}{\begin{example}}
\newcommand{\eeg}{\end{example}}
\newcommand{\begs}{\begin{examples}}
\newcommand{\eegs}{\end{examples}}
\newcommand{\bdefe}{\begin{defn}}
\newcommand{\edefe}{\end{defn}}
\newcommand{\bprob}{\begin{prob}}
\newcommand{\eprob}{\end{prob}}
\newcommand{\bques}{\begin{ques}}
\newcommand{\eques}{\end{ques}}
\newcommand{\bei}{\begin{itemize}}
\newcommand{\eei}{\end{itemize}}

\newcommand{\bde}{\begin{deter}}
\newcommand{\ede}{\end{deter}}
\newcommand{\bca}{\begin{case}}
\newcommand{\eca}{\end{case}}
\newcommand{\bcl}{\begin{claim}}
\newcommand{\ecl}{\end{claim}}
\newcommand{\bcon}{\begin{conj}}
\newcommand{\econ}{\end{conj}}
\newcommand{\bcons}{\begin{conjs}}
\newcommand{\econs}{\end{conjs}}
\newcommand{\bprop}{\begin{propo}}
\newcommand{\eprop}{\end{propo}}
\newcommand{\br}{\begin{rem}}
\newcommand{\er}{\end{rem}}
\newcommand{\brs}{\begin{rems}}
\newcommand{\ers}{\end{rems}}
\newcommand{\bo}{\begin{obser}}
\newcommand{\eo}{\end{obser}}
\newcommand{\bos}{\begin{obsers}}
\newcommand{\eos}{\end{obsers}}
\newcommand{\bpf}{\begin{pf}}
\newcommand{\epf}{\end{pf}}
\newcommand{\ba}{\begin{array}}
\newcommand{\ea}{\end{array}}
\newcommand{\beq}{\begin{eqnarray}}
\newcommand{\beqq}{\begin{eqnarray*}}
\newcommand{\eeq}{\end{eqnarray}}
\newcommand{\eeqq}{\end{eqnarray*}}

\newcommand{\ds}{\displaystyle}

\newcounter{minutes}
\divide\time by 60
\newcounter{hours}
\multiply\time by 60 \addtocounter{minutes}{-\time}

\begin{document}
\title[Generalized Zalcman conjecture for convex functions of order $\alpha$]
{Generalized Zalcman conjecture for convex functions of order $\alpha$}

\thanks{
File:~\jobname .tex,
          printed: \number\day-\number\month-\number\year,
          \thehours.\ifnum\theminutes<10{0}\fi\theminutes}

\author[L. Li]{Liulan Li}
\address{L. Li, College of Mathematics and Statistics,
Hengyang Normal University, Hengyang,  Hunan 421002, People's Republic of China.
}
\email{lanlimail2012@sina.cn}

\author[S. Ponnusamy]{Saminathan Ponnusamy
}
\address{S.  Ponnusamy,
Indian Statistical Institute (ISI), Chennai Centre, SETS (Society
for Electronic Transactions and Security), MGR Knowledge City, CIT
Campus, Taramani, Chennai 600 113, India.
}
\email{samy@isichennai.res.in, samy@iitm.ac.in}

\author[J. Qiao]{Jinjing  Qiao $^\dagger $
}
\address{J. Qiao, Department of Mathematics,
Hebei University, Baoding, Hebei 071002, People's Republic of China.}
\email{mathqiao@126.com}


\keywords{Univalent, starlike, convex and close-to-convex functions, extreme points, closed convex hull and subordination,
Zalcman conjecture. \\
${} ^\dagger$ Corresponding author.
}


\begin{abstract}
Let $\mathcal S$ denote the class of all functions of the form $f(z)=z+a_2z^2+a_3z^3+\cdots$ which are analytic and univalent in the
open unit disk $\ID$ and, for $\lambda >0$, let $\Phi_\lambda (n,f)=\lambda a_n^2-a_{2n-1}$ denote the generalized Zalcman coefficient functional.
Zalcman conjectured that if $f\in \mathcal S$, then  $|\Phi_1 (n,f)|\leq (n-1)^2$ for $n\ge 3$. 
The functional of the form $\Phi_\lambda (n,f)$ is indeed related to Fekete-Szeg\H{o} functional of the $n$-th root transform of the
corresponding function in $\mathcal S$.
This conjecture has been verified for a certain special geometric subclasses of $\mathcal S$ but the conjecture remains open for
$f\in {\mathcal S}$ and for $n > 6$. In the present paper,  we prove sharp bounds on $|\Phi_\lambda (n,f)|$ for $f\in \mathcal{F}(\alpha )$
and for all $n\geq 3$, in the case that  $\lambda$ is a  positive real parameter, where $ \mathcal{F}(\alpha )$ denotes the
family of all functions $f\in {\mathcal S}$ satisfying the condition
$${\rm Re } \left( 1+\frac{zf''(z)}{f'(z)}\right) > \alpha   ~\mbox{ for $z\in \ID$},
$$
where $-1/2\leq \alpha <1$. Thus, the present article proves the generalized Zalcman conjecture for convex functions of order $\alpha$, $\alpha \in [-1/2,1)$.
\end{abstract}

\maketitle \pagestyle{myheadings}
\markboth{L. Li, S. Ponnusamy, and J. Qiao}{Generalized Zalcman conjectures}

\section{Introduction and preliminaries}

Let ${\mathcal H}$ be the set of all analytic functions in the unit disk $\mathbb{D}=\{z\in\IC:\,|z|<1\}$ and
${\mathcal A}=\{f\in {\mathcal H}:\, f(0)=f'(0)-1=0\}$. Clearly each $f\in \mathcal A$ has the form
\be\label{eq1}
f(z)=z+\sum_{n=2}^{\infty}a_{n}z^{n}, \quad z\in\ID.
\ee
For a constant $\alpha\in [-1/2,1)$, a function $f\in \mathcal A$ is said to be in the class $ \mathcal{F}(\alpha )$ if  $f$ satisfies
the condition
\be\label{eq1-e1}
{\rm Re } \left( 1+\frac{zf''(z)}{f'(z)}\right) >\alpha  ~\mbox{ for $z\in \ID$}.
\ee
A number of important properties of the family $ \mathcal{F}(\alpha
)$ for various special values of $\alpha$ may be obtained from the
literature \cite{Du,Go,HM,Pom}. For example,  the family
$\mathcal{F}(0)$ consists of normalized convex functions, usually
denoted by the symbol $\mathcal K$, and thus, for $\alpha\in [0,1)$,
functions in $\mathcal{F}(\alpha )$ are convex in $\ID$. Moreover, a
function $f\in \mathcal A$ is convex  precisely when the function
$g(z)=zf'(z)$ is starlike, i.e. $g(\ID)$ is a domain which is
starlike (with respect to the origin). Thus, $g\in \mathcal A$ is
starlike if ${\rm Re } (zg'(z)/g(z)) >0$ in $\ID$. Also, it is worth
to recall that functions in $\mathcal{F}(-1/2)$  (and hence in $
\mathcal{F}(\alpha )$ for $\alpha\in [-1/2,0)$) are known to be
convex in one direction (and hence functions in $\mathcal{F}(-1/2)$ are close-to-convex and univalent in
$\ID$) but are not necessarily starlike in $\ID$ (see
\cite{Umezawa}). Moreover, if $f\in \mathcal{F}(\alpha )$ and is of
the form \eqref{eq1}, then one has the following necessary
coefficient inequality  (see for instance, \cite[Theorem 5.6,
p.324]{Suf-05}) \be\label{eq1-e2} |a_n|\leq A_n:=\frac{\Gamma
(n+1-2\alpha )}{n!\Gamma(2-2\alpha)} \ee for $\alpha\in [-1/2,1)$
and for all $n\geq 2$, where $\Gamma (\cdot )$ denotes the usual
gamma function.  See the relation \eqref{eq5} in the proof of Lemma
\ref{basic estimate} for a quick proof of the coefficient inequality
\eqref{eq1-e2}. Throughout $A_n:=A_n(\alpha)$ denotes the Taylor
coefficients of the extremal function $f_{\alpha}\in
\mathcal{F}(\alpha )$, where \be\label{eq1-e4}
f_{\alpha}(z)=\frac{1-(1-z)^{2\alpha-1}}{2\alpha-1} \ee and for
$\alpha =1/2$ this is interpreted as the limiting case which gives
$f_{1/2}(z)=-\log (1-z)$. We refer to the recent articles
\cite{BhaPo-JRMS,samy-hiroshi-swadesh} for certain properties of \textit{sections/partial sums} of
functions from the class $\mathcal{F}(-1/2)$. The importance of the
class $\mathcal{F}(-1/2)$ connected with certain univalent harmonic
mappings are considered in \cite{Bshouty-Lyzzaik-2010}.

The aim of this article is to solve the generalized Zalcman coefficient conjecture for the class
$\mathcal{F}(\alpha )$. We begin to present necessary preliminaries.

One of the classical problems is to find for each $\lambda >0$ the maximum modulus value of the
\textit{generalized Zalcman functional}
$$\Phi_\lambda (n,f):=\lambda a_n^2-a_{2n-1}
$$
over the class $\mathcal S$ of functions $f$ of the form \eqref{eq1}. First we remark that the functional for the case 
$n=2$ is fundamental in the investigation of a number of problems in function theory and is popularly known as Fekete and Szeg\H{o} 
functional \cite{FekSze33}. Secondly, we observe that for $\theta \in \IR$,
$$\Phi_\lambda (n,e^{-i\theta}f(e^{i\theta}z))=e^{2(n-1)i\theta}\Phi_\lambda (n,f(z))
$$
and thus, $|\Phi_\lambda (n,f)|$ is invariant under rotation.
For the special case $\lambda =1$,  $\Phi_1 (n,f)$  is simply referred to as the \textit{Zalcman functional}
for $f\in {\mathcal S}$. In 1960, Lawrence Zalcman conjectured that the sharp inequality
\be\label{eq1-e3}
|\Phi_1(n,f)|=|a_n^2-a_{2n-1}|\leq (n-1)^2
\ee
holds for $f\in {\mathcal S}$ and for all $n\geq 3$, with equality
only for the Koebe function $k(z)=z/(1-z)^2$ and its rotations.
This remarkable conjecture, also called the Zalcman coefficient inequality, which was posed as an approach to prove
the Bieberbach conjecture,  was investigated by many mathematicians,
and remained open for all $n > 6$.

By means of Loewners's method,
Fekete and Szeg\H{o} \cite{FekSze33} (see also \cite[Theorem 3.8]{Du}) indeed obtained
the following result. Later in 1960, the same was derived by Jenkins \cite{Jen60} by means of his general coefficient
theorem.

\begin{Thm} $($\cite{FekSze33}$)$ \label{Th1}
For each $f\in {\mathcal S}$, $|\lambda
a_2^2-a_{3}|\leq 1 +2\exp(-2\lambda (1-\lambda))$ for $0\leq \lambda
<1$. The bound is sharp for each  $\lambda$.
\end{Thm}

In particular, for $\lambda \rightarrow 1^{-}$, we have  the well-known Fekete-Szeg\H{o} inequality
(i.e. the case $n=2$ of \eqref{eq1-e3}):
$$|\Phi _1(2,f)|=|a_2^2-a_{3}|\leq 1
$$
(see also \cite[Theorem 1.5]{Pom}).
In 1985, Pfluger \cite{Pflu-85} employed the variational method to present another treatment of Fekete-Szeg\H{o} inequality
and in 1986, Pfluger \cite{Pflu-86} used the method of Jenkins to obtain Theorem~\Ref{Th1} for certain complex values of $\lambda$.
Yet another important remark is that the functional
$\Phi _1(2,f)=a_2^2-a_{3}$ becomes $S_f(0)/6$, where $S_f$ denotes the Schwarzian derivative which is defined for locally univalent function $f$ by
$$S_f=(f''/f')'-(1/2)(f''/f')^2.
$$
Next, if we consider the $n$-th root transform
$$g(z)=\sqrt[n]{f(z^n)} :=z\sqrt[n]{f(z^n)/z^n}=z+c_{n+1}z^{n+1} +c_{2n+1}z^{2n+1} +\cdots
$$
of $f\in {\mathcal S}$ with the power series of \eqref{eq1}, we find that
$$\lambda a_2^2-a_{3}=n(\mu c_{n+1}^2-c_{2n+1}),
$$
where $ \mu =\lambda n +(n-1)/2$. This observation clearly defines the role of the generalized Zalcman functional in the class $\mathcal S$ through the 
Fekete-Szeg\H{o} functional. Thus, generalized Zalcman functional can be regarded as the generalization of Fekete-Szeg\H{o} functional defined 
as in Theorem \Ref{Th1}.

Sharp bound for the generalized Fekete-Szeg\H{o} functional has been established for several subclasses of ${\mathcal S}$ (see
\cite{BT,MaWa,MaWi}) and more recently in \cite{AbLiPo-2014,LiPo-2014,LiPo-PP2016}.  The Zalcman coefficient inequality for $n=3$ and for the full class ${\mathcal S}$, was established in
\cite{Kr95} and also for the special cases $n=4,5,6$ in \cite{Kr10}.  Recently, the authors in \cite{AbLiPo-2014} considered the
generalized Fekete-Szeg\H{o} inequality for $ \mathcal{F}(\alpha)$ and the generalized Zalcman coefficient inequality for the class $ \mathcal{F}(-1/2)$.
We refer to Theorems 2.1, 2.2 and 3.3 in \cite {LiPo-2014} for the precise formulation of these results.

In this note we solve the generalized Zalcman coefficient inequality for the class $ \mathcal{F}(\alpha)$ and
obtain certain earlier known results as corollaries to it (see, for example, Corollary \ref{n3}). In Section \ref{sec4}, we present a number
of lemmas and the main results are stated and proved in Section \ref{sec4a}.

\section{Representation of functions in $ {\mathcal F}(\alpha )$}\label{sec4}

\begin{lem}\label{basic estimate}
Let $-\frac{1}{2}\leq\alpha <1$ and $f\in {\mathcal F}(\alpha )$ as in the form \eqref{eq1}. Then we have
\be\label{eq3}
f'(z)=\int^{2\pi}_{0} (1-e^{i\theta} z)^{2\alpha-2}\,d\nu(\theta),
\ee
where $\nu(\theta)$ is a probability measure on $[0, 2\pi]$ and for $n\geq2$,
\be\label{eq4}
|\lambda a^2_n-a_{2n-1}|\leq  (\lambda A_n^2 - 2A_{2n-1})\int^{2\pi}_{0} \cos^2(n-1)\theta \,d\nu(\theta) + A_{2n-1},
\ee
where $A_n:=A_n(\alpha)$ is given by \eqref{eq1-e2}.
\end{lem}
\bpf
Let $f\in {\mathcal F}(\alpha )$ for some $\alpha \in [-1/2,1)$.  By the well-known Herglotz representation theorem for analytic functions $p$ with positive
real part in $\ID$, $p(0)=1$, and the analytic characterization of $f\in {\mathcal F}(\alpha )$ given by \eqref{eq1-e1}, one has
$$\frac{1}{1-\alpha}\left (1+\frac{zf''(z)}{f'(z)}-\alpha\right )=p(z):=\int^{2\pi}_{0} \frac{1+e^{i\theta} z}{1-e^{i\theta} z}\,d\nu (\theta)  , \quad |z|<1,
$$
where $\nu(\theta)$ is a probability measure on $[0, 2\pi]$.
By a computation, we easily have
\beqq
f'(z)&=&\int^{2\pi}_{0} (1-e^{i\theta} z)^{2\alpha-2}\, d\nu(\theta)\\
&=& 1+\sum^\infty_{n=2}(-1)^{n-1}\frac{(2\alpha-2)(2\alpha-3)\cdots (2\alpha-n)}{(n-1)!}
\left(\int^{2\pi}_{0} e^{i(n-1)\theta}\, d\nu(\theta)\right)z^{n-1}.
\eeqq
By comparing the coefficients of $z^{n-1}$ on both sides of the above equation, we easily have
\be\label{eq5}
a_n=A_n\int^{2\pi}_{0} e^{i(n-1)\theta}\, d\nu(\theta)
~\mbox{ for }~n=2, 3,\ldots,
\ee
where $A_n=A_n(\alpha)$ is given by \eqref{eq1-e2}, i.e.
$$A_n=\frac{\Gamma (n+1-2\alpha )}{n!\Gamma(2-2\alpha)}, \quad n\ge 2.
$$
We observe that the last relation quickly gives the necessary coefficient inequality \eqref{eq1-e2} for functions in ${\mathcal F}(\alpha )$.
Since $|\lambda a^2_n-a_{2n-1}|$ is invariant under rotations, we
can consider instead the problem of maximizing the functional ${\rm Re }\,( \lambda a^2_n-a_{2n-1})$.
Consequently, by \eqref{eq5}, we begin to observe that
\beqq {\rm Re }\left(\lambda a^2_n-a_{2n-1}\right)&=& \lambda A_n^2 \left(\int^{2\pi}_{0}
\cos(n-1)\theta \, d\nu(\theta)\right)^2 -\lambda A_n ^2 \left(\int^{2\pi}_{0} \sin(n-1)\theta
\, d\nu (\theta)\right)^2\\
&& - A_{2n-1}\int^{2\pi}_{0} \cos2(n-1)\theta \, d\nu(\theta).
\eeqq
If we apply Cauchy-Schwarz inequality to the first integral above and use the trigonometric identity
$\cos 2t =2\cos ^2 t -1$ in the third integral, we find that
\beq
{\rm Re }\left(\lambda a^2_n-a_{2n-1}\right)
&\leq&\lambda A_n^2 \int^{2\pi}_{0} \cos^2(n-1)\theta \, d\nu(\theta) -2A_{2n-1}\int^{2\pi}_{0} \cos^2(n-1)\theta \, d\nu(\theta) +A_{2n-1}
\nonumber  \\
\nonumber &=& (\lambda A_n^2 - 2A_{2n-1})\int^{2\pi}_{0} \cos^2(n-1)\theta \, d\nu(\theta)   +A_{2n-1} .
\eeq
The proof is completed.
 \epf


Here is an alternate approach to Lemma \ref{basic estimate} which works for a more general setting.
If $X$ is a linear topological space, then a subset $Y$ of $X$ is called convex if $tx+(1-t)y\in Y$ whenever $x,\ y\in Y$
and $0\leq t\leq 1$. The closed convex hull of $Y$ is defined as the intersection of all closed convex sets containing $Y$. A point $u\in
Y$ is called an extremal point of $Y$ if $u=tx+(1-t)y$, $0<t<1$ and $x, y\in Y$, implies that $x=y$.
See \cite{HM, MacGreW} for a general reference and for many important results on this topic.

In order to solve the generalized Zalcman coefficient inequality problem for the class ${\mathcal F}(\alpha )$, we need the following lemma.

\begin{Lem}\label{extreme} $($\cite{HM}$)$
Suppose that $F_{\alpha}(z)=(1-z)^{2\alpha-2}$ and
$\alpha \in [-1/2,1/2]$. If $s(F_{\alpha})$,
$\mathcal{H}s(F_{\alpha})$ and $\mathcal{EH}s(F_{\alpha})$ denote
the set of analytic functions subordinate to $F_{\alpha}$, the
closed convex hull of $s(F_{\alpha})$ and the set of the extremal
points of $\mathcal{H}s(F_{\alpha})$, respectively, then
$\mathcal{H}s(F_{\alpha})$ consists of all analytic functions
represented by \be\label{eq6} F_{\alpha}(z)=\int_{|x|=1} (1-x
z)^{2\alpha-2}d\mu(x), \ee where $\mu(x)$ is a probability measure
on the unit circle $\partial \mathbb{D}$. Moreover,
$\mathcal{EH}s(F_{\alpha})$ consists of the functions given by
\be\label{eq7} F_{\alpha}(z)=(1-x z)^{2\alpha-2}, \ee where $|x|=1$.
\end{Lem}

If $\mathcal{F}\subset\mathcal{H}$ is convex and $L:\,\mathcal{H}\rightarrow \mathbb{R}$ is a real-valued functional on
$\mathcal{A}$, then we say that $L$ is convex on $\mathcal{F}$ provided that
$$L\left(tg_1+(1-t)g_2\right)\leq t L(g_1)+(1-t)L(g_2)
$$
whenever $g_1,\ g_2\in\mathcal{F}$ and $0\leq t\leq1.$

Since $\mathcal{H}s(F_{\alpha})$ is convex, we have a real-valued, continuous and convex functional on $\mathcal{H}s(F_{\alpha})$.

\begin{lem}\label{convex}
Suppose that $g(z)=1+\sum_{n=2}^{\infty}b_{n}z^{n-1}$ is analytic in $\mathbb{D}$ and
$$J(g)=\lambda\frac{\left({\rm Re\, } b_{n}\right)^2}{n^2}-\frac{{\rm Re\, } b_{2n-1}}{2n-1},
$$
where $\lambda>0$. Then $J$ is a real-valued, continuous and convex
functional on $\mathcal{H}s(F_{\alpha})$.
\end{lem}
\bpf
Let $h(z)=1+\sum_{n=2}^{\infty}c_{n}z^{n-1}$ be analytic in $\mathbb{D}$ and $0\leq t\leq1.$
By the definition of $J$, we have
$$J(h)=\lambda\frac{\left({\rm Re\,}c_{n}\right)^2}{n^2}-\frac{{\rm Re\, } c_{2n-1}}{2n-1}
$$
and thus,
\beqq
J\left(t g+(1-t)h\right)&=&\lambda\frac{t^2\left({\rm Re\, }
b_{n}\right)^2+2t(1-t){\rm Re\, } b_{n}\, {\rm Re\, }
c_{n}+(1-t)^2\left({\rm Re\, }
c_{n}\right)^2}{n^2}\\
& &-\left (\frac{t{\rm Re\, } b_{2n-1}+(1-t){\rm Re\, } c_{2n-1}}{2n-1}\right ).
\eeqq
Therefore, by rearrangements, we find that
$$J\left(t g+(1-t) h\right)-tJ(g)-(1-t)J(h)=-\frac{\lambda (t(1-t)}{n^2}\left({\rm Re\, } b_{n}-{\rm Re\, } c_{n}\right)^2\leq0,
$$
which implies that $J$ is a real-valued, continuous and convex functional on $\mathcal{H}s(F_{\alpha})$.
\epf

Since $s(F_{\alpha})$ is compact, for $J$ as defined in Lemma
\ref{convex}, Theorem 4.6 in \cite{HM} yields the following lemma.

\begin{lem}\label{maximum}
We have
$$\max\{J(f):\, f\in\mathcal{H}s(F_{\alpha})\}=\max\{J(f):\, f\in
s(F_{\alpha})\}=\max\{J(f):\, f\in\mathcal{EH}s(F_{\alpha})\}.
$$
\end{lem}

\begin{lem}\label{upper bound}
Let $-\frac{1}{2}\leq\alpha<\frac{1}{2}$, $\lambda>0$ and $f\in {\mathcal F} (\alpha)$ be as in the form \eqref{eq1}. Then  $f$ has the form
\eqref{eq3} and for $n\geq 2$, we have
\be\label{eq4-e4}
|\lambda a^2_n-a_{2n-1}|\leq (\lambda A_n^2 - 2A_{2n-1})\cos^2(n-1)\theta    +A_{2n-1},
\ee
where $A_n$ is given by \eqref{eq1-e2}.
\end{lem}
\bpf By Lemma \ref{basic estimate},  there exists a function $g$ analytic in $\ID$ such that
$$g(z)=\int^{2\pi}_{0} (1-e^{i\theta} z)^{2\alpha-2}d\nu(\theta)=1+\sum_{n=2}^{\infty}b_{n}z^{n-1},
$$
and $f'(z)=g(z)$ so that $na_n=b_{n}$. Now, Lemma \Ref{extreme} shows that $f'=g\in \mathcal{H}s(F_{\alpha})$.

Since $|\lambda a^2_n-a_{2n-1}|$ is invariant under rotations, we can consider instead the problem of maximizing the functional
${\rm Re }\,( \lambda a^2_n-a_{2n-1})$ so that
\beqq
{\rm Re }\,( \lambda a^2_n-a_{2n-1})
&= &\lambda ({\rm Re }\, a_n)^2-\lambda ({\rm Im }\, a_n)^2-{\rm Re }\, a_{2n-1}\\
&\leq& \lambda ({\rm Re }\, a_n)^2-{\rm Re }\, a_{2n-1}\\
&=&\lambda\frac{\left({\rm Re\, } b_{n}\right)^2}{n^2}-\frac{{\rm Re\, } b_{2n-1}}{2n-1}=J(g).
\eeqq
The above facts and Lemma \ref{maximum} imply that
\be\label{eq8}
|\lambda a^2_n-a_{2n-1}|\leq \max\{J(h):\, h\in\mathcal{H}s(F_{\alpha})\}=\max\{J(h):\,
h\in\mathcal{EH}s(F_{\alpha})\}.
\ee
Since each $h\in\mathcal{EH}s(F_{\alpha})$ is in the form
$$h(z)=(1-e^{i\theta} z)^{2\alpha-2}=1+\sum^\infty_{n=2}nA_ne^{i(n-1)\theta}z^{n-1},
$$
as in the proof of Lemma \ref{basic estimate}, the relation \eqref{eq8} reduces to
$$|\lambda a^2_n-a_{2n-1}| \leq (\lambda A_n^2 - 2A_{2n-1})\cos^2(n-1)\theta  +A_{2n-1}
$$
and the proof is complete.
\epf


\begin{lem}\label{monotonical}
For $-\frac{1}{2}\leq\alpha<1$ and $n\geq3$, we define
\be\label{eq1-e6}
C_n(\alpha)=\frac{2A_{2n-1} (\alpha )}{A_n^2(\alpha )},
\ee
where $A_n=A_n(\alpha )$ is given by \eqref{eq1-e2}. Then for fixed $\alpha$, we have
\bee
\item [{\rm (1)}] $C_n(\alpha)$ is monotonically decreasing with respect to $n$ and $C_n(\alpha)\leq C_3(\alpha)$ if $-\frac{1}{2}\leq\alpha<0$;
\item [{\rm (2)}] $C_n(\alpha)$ is monotonically increasing with respect to $n$ and $C_n(\alpha)\geq C_3(\alpha)$ if $0<\alpha<1$,
\eee
where
$$C_3(\alpha)=\frac{3}{5}\frac{(2\alpha-4)(2\alpha-5)}{(2\alpha-2)(2\alpha-3)}.
$$
\end{lem}
\bpf  Note that $C_n(0)=2$ for all $n\ge 3$ and thus, there is nothing to prove for the case $\alpha =0$.
Clearly, $C_n(\alpha)>0$ for $-\frac{1}{2}\leq\alpha<1$. We only need to consider
\beqq
\frac{C_{n+1}(\alpha)}{C_n(\alpha)}-1&=&\frac{(n+1)^2(2n+1-2\alpha)(n-\alpha)}{n(2n+1)(n+1-2\alpha)^2}-1\\
&=&\frac{\alpha [(n+1)(4n^2-n-1)-2\alpha(3n^2-1)]}{n(2n+1)(n+1-2\alpha)^2}\\
&=&\frac{\alpha \varphi (n) }{n(2n+1)(n+1-2\alpha)^2},
\eeqq
where $\varphi (n)=4n^3+3n^2(1-2\alpha)-2n-1+2\alpha$. Since  $\varphi (n)$ is trivially an increasing function of $n$ for $n\ge 3$
and for each $-\frac{1}{2}\leq\alpha<1$, it follows that $\varphi (n)\geq \varphi (3)=128-52\alpha >0$. This observation shows that
$C_{n+1}(\alpha)> C_n(\alpha)$ for $0<\alpha<1$ while $C_{n+1}(\alpha)< C_n(\alpha)$ if $-\frac{1}{2}\leq\alpha<0$. The desired conclusion now follows.
\epf

\section{Main results}\label{sec4a}

Let $Co\,({\mathcal F})$ denote the convex hull of the set ${\mathcal F}$ and its closure by $\overline{Co}\,({\mathcal F})$.
In view of the extreme points method described in Lemma \ref{upper bound}, our results continue to hold if we replace
the assumption $f\in {\mathcal F}(\alpha)$ by $f\in \overline{Co}\,({\mathcal F}(\alpha ))$ for the case $-\frac{1}{2}\leq\alpha <\frac{1}{2}$.
For example, we have the following results and for the sake of completeness we include the proofs here.

\begin{thm}\label{0}
For ${\mathcal F} (0):={\mathcal K}$, let $f\in \overline{Co}\,({\mathcal K})$ as in the form \eqref{eq1}.
Then we have
\bee
\item [{\rm (1)}] $|\lambda a_n^2-a_{2n-1}|\leq \lambda-1$ for $n\geq 3$ and $\lambda\geq2$. The
equality is attained for the function $ f_0(z)=\frac{z}{1-z}$.
\item[{\rm (2)}] $|\lambda a_n^2-a_{2n-1}|\leq 1$ for $n\geq 3$ and $0<\lambda<2$. The
equality is attained by convex combination of convex functions in
$\mathcal{K}$, namely, for the functions $f$ in the form
$$f(z)=\sum^{2n-3}_{k=0}\alpha_k \frac{z}{1-e^{i\theta_k}z}, 
$$
where $0\leq \alpha_k\leq1$,
$$\theta_k=\frac{(2k+1)\pi}{2n-2} ~\mbox{ and }~ \sum^{n-2}_{m=0}\alpha_{2m}=\sum^{n-2}_{m=0}\alpha_{2m+1}=\frac{1}{2}.
$$
\eee
\end{thm}
\bpf For $\alpha=0$, we have $A_n=1$ for all $n\geq 2$. Thus we can apply either Lemma \ref{basic estimate} or Lemma \ref{upper bound}.
In either way, applying either \eqref{eq4} or \eqref{eq4-e4} we see that for $\lambda\geq2$, we have
$$|\lambda a^2_n-a_{2n-1}|\leq (\lambda -2)\cdot 1 +1= \lambda -1,
$$
where the equality is attained by the convex function $f_0(z)=\frac{z}{1-z}$. For $0<\lambda<2$, both \eqref{eq4} and \eqref{eq4-e4} reduces to
$$|\lambda a^2_n-a_{2n-1}|\leq 1,
$$
which occurs when $\theta=\pi /(2(n-1))$ and at this point $\sin ^{2}((n-1)\theta)=1$.
\epf

Theorem \ref{0}(2) is also obtained recently in \cite{EfrVuko-PP2014} (see also \cite{AgrSahoo-PP}). Next we consider the case $\alpha =-1/2$
and because of its independent interest we supply the proof.

\begin{cor}\label{n3} $($\cite[Theorem 3.3]{LiPo-2014}$)$
 Let $f\in \mathcal{F}(-1/2)$ and $f$ be of the
form \eqref{eq1}. Then we have \bee
\item [{\rm (1)}] $|\lambda a_n^2-a_{2n-1}|\leq\frac{ (n+1)^2}{4}\lambda-n$ for $n\geq 3$ and $\lambda\geq\frac{3}{2}$. The
equality is attained for the function $ f_{-1/2}(z)$ given by
\be\label{eqf0}
f_{-1/2}(z)=\frac{z-z^2/2}{(1-z)^2} = z+\sum_{n=2}^\infty \frac{1+n}{2}z^n.
\ee
\item[{\rm (2)}] $|\lambda a_n^2-a_{2n-1}|\leq\frac{ (n+1)^2}{4}\lambda-n$ for $0<\lambda<\frac{3}{2}$ and
$n>\frac{4-\lambda+2\sqrt{4-2\lambda}}{\lambda}$. The equality is
attained for the function $ f_{-1/2}(z)$ given by \eqref{eqf0}.

\item[{\rm (3)}] $|\lambda a_n^2-a_{2n-1}|\leq n$ for $0<\lambda<\frac{3}{2}$ and
$3\leq n\leq\frac{4-\lambda+2\sqrt{4-2\lambda}}{\lambda}$. The
equality is attained for  functions $f$ in the following form
$$f(z)=\sum^{2n-3}_{k=0}\alpha_k \frac{2z-e^{i\theta_k}z^2}{2(1-e^{i\theta_k}z)^2}=\sum^{2n-3}_{k=0}\alpha_k e^{-i\theta_k} f_{-1/2}(ze^{i\theta_k}) ,
$$
where $0\leq \alpha_k\leq1$,
$$\theta_k=\frac{(2k+1)\pi}{2n-2} ~\mbox{ and }~ \sum^{n-2}_{m=0}\alpha_{2m}=\sum^{n-2}_{m=0}\alpha_{2m+1}=\frac{1}{2}.
$$
\eee
\end{cor}
\bpf
Set $\alpha =-1/2$ in Lemma \ref{basic estimate} or  Lemma \ref{upper bound}, or apply Lemma \ref{monotonical} directly.
Then, because $A_n(-1/2) =(n+1)/2$ for all $n\ge 2$,  it is clear from \eqref{eq4} or \eqref{eq4-e4} that
$$|\lambda a^2_n-a_{2n-1}|\leq \left \{
\begin{array}{rl}  \ds \left (\lambda \frac{(n+1)^2}{4} -2n\right )+n & \mbox{ if $\lambda\geq \frac{8n}{(n+1)^2}$}\\
 \ds n & \mbox{ if $0<\lambda\leq  \frac{8n}{(n+1)^2}$.}
\end{array} \right .
$$
Clearly, $\frac{8n}{(n+1)^2}\leq \frac{3}{2}$ if and only if $(3n-1)(n-3)\ge 0$ and thus, the Case (1) follows. For  $0<\lambda<\frac{3}{2}$, we see that
$\lambda \ge \frac{8n}{(n+1)^2}$ if and only if $\varphi (n):=\lambda n^2-2n(4-\lambda )+\lambda \geq 0$. Since
$$\varphi (n)= \lambda \left [n-\left (\frac{4-\lambda+2\sqrt{4-2\lambda}}{\lambda}\right )
\right ] \left [n-\left (\frac{4-\lambda-2\sqrt{4-2\lambda}}{\lambda}\right )\right ],
$$
Case (2) follows. Finally, in the last case the range of $n$ shows that $\lambda < \frac{8n}{(n+1)^2}$ and thus, Case (3) is clear.
\epf

Setting $\lambda=1$ and $\frac{3}{2}$ in Corollary \ref{n3}, we obtain the following.

\begin{cor}
Let $f\in \mathcal{F}(-1/2)$ as in the form \eqref{eq1}. Then we have
\bee
\item[{\rm (1)}] $|a_n^2-a_{2n-1}|\leq\frac{ (n-1)^2}{4}$ for $n>5$. The equality is
attained for the function $ f_{-1/2}(z)$ given by \eqref{eqf0}.

\item[{\rm (2)}] $|a_n^2-a_{2n-1}|\leq n$ for $3\leq n\leq5$.

\item[{\rm (3)}] $\ds \big |\frac{3}{2}a_n^2-a_{2n-1}\big |\leq\frac{3n^2-2n+3}{8}$ for $n>3$. The equality is
attained for the function $ f_{-1/2}(z)$ given by \eqref{eqf0}.

\item[{\rm (4)}] $\ds \big |\frac{3}{2}a_3^2-a_5\big |\leq 3$.
\eee
\end{cor}

\begin{thm}\label{less than 0}
Let $n\geq 3$, $-\frac{1}{2}\leq\alpha <0$, $f\in {\mathcal F}(\alpha)$ as in the form \eqref{eq1},
$A_n=A_n(\alpha )$ and $C_n=C_n(\alpha )$  be given by \eqref{eq1-e2} and  \eqref{eq1-e6}, respectively.
\bee
\item [{\rm (1)}] If $n\geq 3$ and $\lambda\geq C_3(\alpha)$, then we have
$$|\lambda a^2_n-a_{2n-1}|\leq  \lambda A_n^2 - A_{2n-1},
$$
where the equality is attained for the function $f_{\alpha}(z)$ defined by \eqref{eq1-e4}.

\item [{\rm (2)}]  If $0<\lambda<C_3(\alpha)$, then there exists a fixed
$n_0>3$ such that
$$C_{n_0-1}(\alpha)>\lambda\geq C_{n_0}(\alpha).
$$
If $0<\lambda<C_3(\alpha)$ and $n\geq n_0$, then
$$|\lambda a^2_n-a_{2n-1}|\leq \lambda A_n^2 - A_{2n-1},
$$
where the equality is attained for the function $f_{\alpha}(z)$ defined by \eqref{eq1-e4}.

\item [{\rm (3)}] If $0<\lambda<C_3(\alpha)$ and $3\leq n<n_0$, then
$$|\lambda a^2_n-a_{2n-1}|\leq A_{2n-1},
$$
where the equality is attained by convex combination of rotations of functions $f_{\alpha}\in{\mathcal F}(\alpha)$.
\eee
\end{thm}
\bpf We apply  Lemma \ref{basic estimate} or Lemma \ref{upper bound}. Thus, by \eqref{eq4} or \eqref{eq4-e4}, we find that
\be\label{eq1-e5}
|\lambda a^2_n-a_{2n-1}|\leq \left \{\begin{array}{rl}  \ds \lambda A_n^2 - A_{2n-1} & \mbox{ if $\lambda A_n^2 - 2A_{2n-1}\geq 0$}\\
 \ds A_{2n-1}  & \mbox{ if $\lambda A_n^2 - 2A_{2n-1}\leq 0$.}
\end{array} \right .
\ee
This is the key and using this and Lemma \ref{monotonical}, we obtain the desired conclusion in each case.
We remind that the inequality \eqref{eq1-e5} holds for $-\frac{1}{2}\leq\alpha <1$.\\

\noindent {\bf Case (1)}. $\lambda\geq C_3(\alpha)$ and
$n\geq3$.\smallskip

In this case, $\lambda\geq C_3(\alpha)\geq C_n(\alpha)$ by Lemma \ref{monotonical} and thus, $\lambda A_n^2 - 2A_{2n-1}\geq 0$ for all $n\ge 3$, which implies
the conclusion of Case (1), where the equality is attained by the function $f_\alpha(z)$.

If $0<\lambda<C_3(\alpha)$, then there exists a fixed $n_0>3$ such that
$$C_{n_0-1}(\alpha)>\lambda\geq C_{n_0}(\alpha),
$$
by Lemma \ref{monotonical}. So we need to divide the case $0<\lambda<C_3(\alpha)$ into the following two cases.\\

\noindent {\bf Case (2)}. $0<\lambda<C_3(\alpha)$ and $n\geq n_0$.
\smallskip

In this case, Lemma \ref{monotonical} yields that
$\lambda A_n^2 - 2A_{2n-1}\geq0$ for  $n\geq n_0$, and the conclusion follows from \eqref{eq1-e5}.\\

\noindent {\bf Case (3)}. $0<\lambda<C_3(\alpha)$ and $3\leq n<n_0$.\smallskip

In this case, $\lambda A_n^2 - 2A_{2n-1}<0$ for $3\leq  n<n_0$ and the desired inequality follows from \eqref{eq1-e5}.
\epf

\vspace{6pt}

Proof of the following theorem is similar and it just uses Lemma \ref{monotonical} and the equation \eqref{eq1-e5}. Thus, we include only the necessary details.

\begin{thm}\label{bigger than 0}
Let $n\geq 3$, $0<\alpha <1$, $\alpha\neq \frac{1}{2}$ and $f\in
{\mathcal F}(\alpha)$ as in the form \eqref{eq1},$A_n=A_n(\alpha )$ and $C_n=C_n(\alpha )$  be given by \eqref{eq1-e2} and  \eqref{eq1-e6}, respectively.

\bee
\item [{\rm (1)}] If $n\geq 3$ and $0<\lambda\leq C_3(\alpha)$, then
$$|\lambda a^2_n-a_{2n-1}|\leq A_{2n-1},
$$
where the equality is attained by
convex combination of rotations of functions $f_{\alpha}\in{\mathcal F}(\alpha)$ defined by \eqref{eq1-e4}.

\item [{\rm (2)}] If $\lambda>C_3(\alpha)$, then there exists a fixed $n_0>3$ such that
$$C_{n_0-1}(\alpha)<\lambda\leq C_{n_0}(\alpha).
$$
Furthermore,  if $\lambda>C_3(\alpha)$ and $n\geq n_0$, then
$$|\lambda a^2_n-a_{2n-1}|\leq A_{2n-1},
$$
where the equality is attained by convex combination of rotations of functions $f_{\alpha}\in{\mathcal F}(\alpha)$ defined by \eqref{eq1-e4}.

\item [{\rm (3)}]  If $\lambda>C_3(\alpha)$ and $3\leq n<n_0$, then
$$|\lambda a^2_n-a_{2n-1}|\leq \lambda A_n^2 - A_{2n-1},
$$
where the equality is attained for the function $f_{\alpha}(z)$ defined by \eqref{eq1-e4}.
\eee
\end{thm}
\bpf  \noindent {\bf Case (1)}. $0<\lambda\leq C_3(\alpha)$ and
$n\geq3$.\smallskip

In this case, $\lambda\leq C_3(\alpha)\leq C_n(\alpha)$ by Lemma
\ref{monotonical} and thus, $\lambda A_n^2 - 2A_{2n-1}\leq 0$, which by \eqref{eq1-e5} implies the desired inequality
$$|\lambda a^2_n-a_{2n-1}|\leq  A_{2n-1}.
$$

If $\lambda>C_3(\alpha)$, then there exists a fixed $n_0>3$ such that $C_{n_0-1}(\alpha)<\lambda\leq C_{n_0}(\alpha)$ by Lemma
\ref{monotonical}. \\

\noindent {\bf Case (2)}. $\lambda>C_3(\alpha)$ and $n\geq n_0$.
\smallskip

For this case, Lemma \ref{monotonical} yields that
$\lambda A_n^2 - 2A_{2n-1}\leq 0$ and the conclusion follows from \eqref{eq1-e5}.\\

\noindent {\bf Case (3)}. $\lambda>C_3(\alpha)$ and $3\leq
n<n_0$.\smallskip

In this case, $\lambda A_n^2 - 2A_{2n-1}\geq0$ and the desired inequality follows as before.
\epf

\begin{thm}
Let $f\in {\mathcal F} (1/2)$ as in the form \eqref{eq1}.
Then we have \bee
\item [{\rm (1)}] $|\lambda a_n^2-a_{2n-1}|\leq \frac{1}{2n-1}$ for $n\geq 3$ and $0<\lambda\leq\frac{18}{5}$,
where the equality is attained by convex combination of rotations of
functions $f_{1/2}\in{\mathcal F}(1/2)$ defined by \eqref{eq1-e4}.

\item[{\rm (2)}] $|\lambda a_n^2-a_{2n-1}|\leq \frac{1}{2n-1}$ for $\lambda>\frac{18}{5}$ and
$n\geq\frac{\lambda+\sqrt{\lambda^2-2\lambda}}{2}$, where the equality is attained by
convex combination of rotations of functions $f_{1/2}\in{\mathcal
F}(1/2)$.

\item[{\rm (3)}] $|\lambda a_n^2-a_{2n-1}|\leq \frac{\lambda}{n^2}-\frac{1}{2n-1}$ for $\lambda>\frac{18}{5}$ and
$3\leq n<\frac{\lambda+\sqrt{\lambda^2-2\lambda}}{2}$. The equality
is attained for the function $f_{1/2}(z)=-\log (1-z)$.
\eee \end{thm}
\bpf By Lemma \ref{basic estimate}, we have $$|\lambda
a^2_n-a_{2n-1}|\leq
\left(\frac{\lambda}{n^2}-\frac{2}{2n-1}\right)\int^{2\pi}_{0}
\cos^2(n-1)\theta \, d\nu(\theta)+\frac{1}{2n-1}.$$

If $n\geq 3$ and $0<\lambda\leq\frac{18}{5}$, then
$$\frac{\lambda}{n^2}-\frac{2}{2n-1}\leq0,
$$
and the desired conclusion holds.

Similarly, if $\lambda>\frac{18}{5}$ and $n\geq\frac{\lambda+\sqrt{\lambda^2-2\lambda}}{2}$, then
$$\frac{\lambda}{n^2}-\frac{2}{2n-1}\leq0,
$$
which shows that the conclusion for this case holds.

If $\lambda>\frac{18}{5}$ and $3\leq n<\frac{\lambda+\sqrt{\lambda^2-2\lambda}}{2}$, then
$$\frac{\lambda}{n^2}-\frac{2}{2n-1}>0,
$$
which provides a proof for this case.
 \epf

\subsection*{Acknowledgements}
This work was completed during the visit of the first author to
Syracuse University. She thanks the university for its hospitality.
The visit and the research of the first author was supported by CSC
of China (No. 201308430274). The research was also supported by NSF
of China (No. 11201130 and No. 11571216), Hunan Provincial Natural
Science Foundation of China (No. 14JJ1012) and construct program of
the key discipline in Hunan province. The second author is on leave
from IIT Madras. The third author
was supported by NSF of China (No. 11501159).

\end{document}